\documentclass[12pt,english]{article}
\usepackage{amsfonts,here,graphicx}
\usepackage{cite,amsfonts,amssymb} 
\usepackage[latin2]{inputenc}
\usepackage[T1]{fontenc}
\def\sup{\mathop{\rm sup}}

\def\proof{\noindent \medskip {\bf Proof:}$\;\;$}

\def\ass#1#2\endass{\vskip5pt plus2pt \noindent{\bf (A.#1)} {\it #2} \vskip5pt
plus2pt }

\newtheorem{lem}{Lemma}

\newtheorem{theor}{Theorem}

\newtheorem{cor}{Corollary}

\newtheorem{defin}{Definition}

\date{} 

\title{\bf Maximal type inequalities for linear stochastic Volterra equations}

\author{\large\sf Anna Karczewska \\
 \\
 Department of Mathematics,
 University of Zielona G\'ora\\
 ul. Szafrana 4a, 65-246 Zielona G\'ora, Poland\\
 e-mail: A.Karczewska@im.uz.zgora.pl\\
}

\begin{document}

\maketitle

\def\thefootnote{}

\vspace{-10mm}
\begin{abstract}
The note is devoted to estimates for convolutions appearing in 
some class of stochastic Volterra equations.
Two maximal inequalities and exponential tail estimate are proved by 
the fractional method of infinite dimensional stochastic calculus. 
The paper extends on non-semigroup case
some results obtained earlier for semigroups.
\end{abstract}

\section{The aim of the paper}\label{aim}

Assume that $(\Omega,\mathcal{F},(\mathcal{F}_t)_{t\geq 0},P)$ is a probability
space with a complete right-continuous filtration and $W(t),~t\geq 0$, 
a cylindrical $(\mathcal{F}_t)$ -- Wiener process with values in a 
separable Hilbert space $U$ and a covariance operator $Q$.
Let $H$ be a separable  Hilbert space with a scalar product $(\cdot,\cdot)$ 
and a norm $|\cdot|$ and let $\{ e_k\}$ be a complete orthonormal system in $H$. 
Assume that $\psi(t),~t\geq 0$, is an appropriate process, defined below. 

In the paper we study $H$-valued stochastic convolutions corresponding to  
linear stochastic Volterra equations of the form 
\begin{equation} \label{e1}
X(t) = X_0 + \int_0^t a(t-\tau)\,AX(\tau)d\tau +
\int_0^t \psi(\tau)\,dW(\tau) \;,
\end{equation}
where $t\in\mathbb{R}_+,~X_0\in H,~a\in L^1_{\mathrm{loc}}
(\mathbb{R}_+),~A$ is a closed linear unbounded operator in $H$ with 
a dense domain $D(A)$ and $W$ is as above.

We assume that the equation (\ref{e1}) is well-posed and denote by 
$\{ S(t) \}_{t\geq 0}\subset B(H)$, where $S(t)(D(A))\subset D(A)$,
the family of bounded linear operators in the space $H$ called 
{\em resolvent} for the equation (\ref{e1}). 
Then the mild solution to  (\ref{e1}) has the form
\begin{equation} \label{e2}
X(t) = S(t)X_0 + \int_0^t S(t-\tau)\,\psi(\tau)\,dW(\tau)\quad\quad t\geq 0  \;.
\end{equation}

The aim of the paper is to provide some estimates for the 
stochastic convolution arising in the mild solution 
(\ref{e2}). In order to do it we will use the {\em factorization method} 
of infinite dimensional stochastic calculus.

Till now some people applied that method for obtaining, among others, the
following results: continuity of mild solutions to stochastic evolution
equations, maximal inequalities or some exponential tail estimates for
stochastic convolutions, see \cite{Za} and references therein. In all papers 
semigroups of operators played the crucial and indispensable role.

In our case, the operators $S(t),~t\geq 0$, do not form any semigroup
and in the consequence, we can not use the known results directly.
Unfortunately, because of the lack of semigroup property, the method used in the
paper does not provide existence of continuous modification of the stochastic 
convolutions considered.
To the best of our knowledge there are no papers joining the factorization
method with stochastic Volterra equations.

\section{Factorization method} \label{fac}

The factorization method in stochastic case consists in
representing trajectories of a process under consideration 
like result of the composition of two
fractional integral operators. When one of them has very smoothing property, 
the whole composition is regular.

The {\em stochastic factorization} method has been 
introduced by DaPrato, Kwapie\'n and Zabczyk \cite{DKZ}. In that method
the crucial role plays a $C_0$-semigroup $R(t),~t\geq 0$, of
bounded, linear operators on a separable Hilbert space $H$ with the
infinitesimal generator $B$, where $D(B)\subset H$. 
For an $H$-valued integrable function
$f,~\alpha\geq 0$, the {\em generalized Riemann--Liouville integral} is defined
as follows
\begin{equation} \label{e3}
 I_\alpha \,f(t) := \frac{1}{\Gamma(\alpha)} \int_0^t 
 R(t-s)(t-s)^{\alpha-1} f(s) ds,\quad t\in[0,T],~T>0  \,.
\end{equation}
The family $I_\alpha,~\alpha>0$, forms a semigroup 
of operators, that is,
$
  I_{\alpha+\beta} = I_\alpha(I_\beta f)(t), \quad t\in [0,T],~~
  \alpha,\beta>0\;.
$

The space $\mathcal{H}_W:=Q^{1/2}U$ is the reproducing kernel of the process $W$.
Let $L_2=L_2(\mathcal{H}_W,H)$ denote the space of all 
Hilbert-Schmidt operators acting from $\mathcal{H}_W$ 
into the space $H$ with the Hilbert-Schmidt norm $||\cdot||_2$.\\[3mm]
($\psi 1$): \hfill \parbox[t]{145mm}{\it
Assume that $\psi(t),\,t\in [0,T]$, is an $L_2$-valued 
predictable process and $\mathbb{E}(\int_0^T ||\psi(t)||_2^2\,dt)<+\infty$.}\\

Then for arbitrary $C_0$-semigroup $R(t),~t\geq 0$,
the stochastic integral
$$
  W_B^\psi(t) := \int_0^t R(t-\tau)\psi(\tau)
  \, dW(\tau),\quad t\in [0,T],
$$
called {\em stochastic convolution}, is a well-defined, $H$-valued stochastic
process. 

In fact, the above stochastic convolution $W_B^\psi$ is well defined under the weaker 
condition on the process $\psi$, that is,
$\mathbb{P}(\int_0^T ||\psi(t)||_2^2\,dt<+\infty)=1$. 
Nevertheless, we introduce the stronger one $(\psi 1)$ because it guarantees
the useful property (see (\ref{e16})) of stochastic integral. 
Then we may write 
$
  W_B^\psi(t) = I_\alpha (I_{1-\alpha}\stackrel{\circ}{W})(t)\;,
$
where $\stackrel{\circ}{W}$ denotes the derivative of the process $W$ and
the process $Y_\alpha:=I_{1-\alpha}\stackrel{\circ}{W}$ 
is defined like the stochastic
integral
$$
  \widetilde{Y}_\alpha (t) = \frac{1}{\Gamma(1-\alpha)} \int_0^t (t-\tau)^{-\alpha}
  R(t-\tau)\psi(\tau)\,dW(\tau), \quad t\in[0,T]\;.
$$
Then the formula
$
  W_B^\psi(t) = I_\alpha \widetilde{Y}_\alpha (t),~~ t\in[0,T],
$ 
is the required stochastic factorization formula.

Now, $W_B^\psi$ is a well-defined $C[0,T]$-valued random variable, provided 
$\widetilde{Y}_\alpha$ is an $L^p(0,T)$-valued random variable with $\alpha>1/p$,
because the operator $I_\alpha$ acts from $L^p(0,T)$ into $C(0,T)$
continuously.

\section{Auxiliary estimates} \label{sc3}

In this section we study relationships between the following processes:
\begin{equation} \label{e4}
 Y(t):= \int_0^t (t-s)^{\beta-1}\, S(t-s)\, \psi(s)\,dW(s) =
  \int_0^t (t-s)^{-\alpha} \,S(t-s) \,\psi(s)\,dW(s) \;,
\end{equation}
\begin{equation} \label{e4a}
 Y_\alpha(t):= \frac{1}{\Gamma(1-\alpha)} \,Y(t)\;,
\end{equation}
\begin{equation} \label{e5}
Z_1(t) := \int_0^t (t-s)^{\alpha-1} \,S(t-s)\, Y(s)\,ds \;,
\end{equation}
\begin{equation} \label{e6}
\mbox{and} \quad\quad Z_2(t) := C_\alpha \int_0^t S(t-s)\,\psi(s)\,dW(s)\;,
\end{equation}
for $t\in [0,T]$, 
where $\alpha,\beta$ are positive numbers such that $\alpha+\beta=1$,
$C_\alpha=\Gamma(\alpha)\Gamma(1-\alpha)=\frac{\pi}{\sin \pi\alpha}$
and $S,\psi,W$ are like earlier with $\psi(s):\; U\rightarrow D(A) \subset H,
\;\; s\geq 0$, satisfying $(\psi 1)$.

In the formulas (\ref{e4})-(\ref{e6}), $S(t),t\geq 0$, denote the corresponding
resolvent operators for Volterra equations of the form (\ref{e1}).
Let us recall that $S(t)$ is linear for each $t\geq 0$, $S(0)\,x=x$ holds on 
$D(A)$, and $S(t)\,x$ is continuous on $\mathbb{R}_+$ for any $x\in D(A)$.
Moreover, $S(t)$ is uniformly bounded on compact intervals. Finally,  $S(t)$
commutes with $A$, that is $S(t)(D(A))\subset D(A)$ and 
$AS(t)x=S(t)Ax$ for all $x\in D(A)$ and $t\geq 0$. Additionally, the so called 
{\em resolvent equation} holds
$ S(t) x = x + \int_0^t a(t-\sigma)AS(\sigma) x d\sigma  $
for all $x\in D(A))$, $t\geq 0$.
By $||S(t)||$ we will denote the norm of the operator $S(t)$, for $t\geq 0$.
For more details concerning such operators we refer to the monograph \cite{Pr}.

For simplicity we assume that the operator $A$ in the equation (\ref{e1}) is
negative and diagonal with respect to the basis $\{ e_k\}$, that is 
$
  A \,e_k = -\mu_k\,e_k,~ \mu_k>0,~ k\in N\;.
$

Let $s(t;\gamma)$ denote the solution of the one-dimensional Volterra equation
\begin{equation} \label{e7}
  s(t;\gamma) + \gamma \int_0^t a(t-\tau)\,s(\tau;\gamma)\,d\tau = 1,\quad
   t\geq 0 \quad \gamma\geq 0 \;,
\end{equation}
where $a \!\in\! L_{\mathrm{loc}}^1(\mathbb{R}_+)$ is the same like in 
(\ref{e1}). Under our assumptions
concerning the resolvent of the equation (\ref{e1}) and the operator $A$, 
the ope\-ra\-tors  $S(t),t\geq 0$, are determined as follows
\begin{equation} \label{e8}
 S(t) \,e_k = s(t;\mu_k)\,e_k, \quad   k\in N\;.
\end{equation}

In the paper we shall study the class of linear Volterra equations of the 
form (\ref{e1}) which satisfy the below hypothesis. 
\\

\noindent{\bf Hypothesis}

\noindent
(s) \hfill\parbox[t]{150mm}{\it The solutions to the equation (\ref{e7}),
 connected with 
the equation (\ref{e1}) are submultiplicative functions, that is for any 
$t,\tau\in [0,T]$, $ s(t+\tau) \leq s(t)\,s(\tau)$.} \\[2mm]
\noindent{\bf Comment:} Integrodifferential equations of the form 
\begin{equation} \label{e8a}
 X(t,\theta) = X_0(\theta)+\frac{1}{\Gamma(\alpha)} \int_0^t (t-s)^{\alpha-1}
\Delta\, X(s,\theta)\,ds\;, 
\end{equation}
where $\Gamma(\alpha)$ is the gamma function, $\Delta$ is Laplacian and 
$\alpha\in [1,2)$, are examples of equations satisfying the above assumption (s).
For more details, see e.g.\ \cite{Pr} or \cite{Fu}.

\begin{cor} \label{c1}
 If $S(t),t\geq 0$, and $s(t;\gamma)$ are like in (\ref{e7})-(\ref{e8}) and the
 assumption (s) holds, then for any $x\in D(A)$
 $$ |S(t+\tau)\,x| \leq |S(t)\,S(\tau)\,x| \;. $$
 Analogously, for any functional $\phi\in H^*$,
\begin{equation} \label{e9}
  \phi (S(t+\tau)\,x) \leq \phi (S(t)\,S(\tau)\,x) \;. 
\end{equation} 
\end{cor}

In order to prove the corollary it is enough to use the relationship (\ref{e8}),
the assumption (s) and linearity of the operators $S(t),t\geq 0$.
\begin{lem} \label{l1}
 Assume that $0<(1/p)<\alpha<1$, the process $Y(t)$ given by (\ref{e4}) is
 well-defined and has $p$-integrable trajectories. When the condition (s) holds,
 then for all $t\in [0,T]$, for any $\phi\in H^*$
\begin{equation} \label{e10}
  \phi (Z_2(t)) \leq \phi (Z_1(t)) \;,
\end{equation} 
where $Z_1(t),Z_2(t),\,t\geq 0$, are defined by formulas  (\ref{e5}) and  (\ref{e6}),
respectively.
\end{lem} 

\proof{
Under the assumption ($\psi 1)$, the stochastic integral 
$\int_0^t \psi(s)\,dW(s)$, for $t\in[0,T]$, is an $H$-valued local martingale.

We introduce the following notation: 
$\phi$ is an arbitrary linear functional belonging to the space $H^*$,
$m(t):= \phi (\int_0^t S(v-s)\psi(s)\,dW(s))$ denotes an auxiliary
square integrable martingale defined for any $t\in [0,v]$, where $v\in [0,T]$.

For any linear functional $\phi$ and $v$ we may write: 
\begin{eqnarray}
\phi(Z_2(v))& = &\Gamma(\alpha)\,\Gamma(1-\alpha)\,m(v) = 
\Gamma(\alpha)\,\Gamma(\beta)\,\int_0^v dm(r) \nonumber\\
 & = & \int_0^v \left[ \int_r^v (v-s)^{\alpha-1} (s-r)^{\beta-1}  ds\right]
 dm(r)  = \nonumber \\
 && \hspace{10ex} \mbox{(from~Fubini's~theorem~for~martingales)} \nonumber\\
 & = & \int_0^v (v-s)^{\alpha-1} \left[ \int_0^s (s-r)^{\beta-1} dm(r)\right]ds
 \nonumber \\
 & = &  \int_0^v (v-s)^{\alpha-1}\left[ \phi \left( \int_0^s (s-r)^{\beta-1} 
 S(v-r)\psi(r)\,dW(r)\right) \right] ds  \nonumber \leq \nonumber\\
 && \hspace{10ex} \mbox{(from~the~property~(\ref{e9}))} \nonumber\\
 & \leq & \int_0^v (v-s)^{\alpha-1}\left[\phi\left( S(v-s) 
 \int_0^s (s-r)^{\beta-1} S (s-r) \psi(r)\,dW(r)\right) \right] ds \nonumber \\
 & = &  \int_0^v (v-s)^{\alpha-1}\phi \left(S(v-s)\,Y(s)\right)ds 
  =\phi(Z_1(v))\;. \nonumber
\end{eqnarray}
\hfill $\blacksquare$  } 

\begin{cor} \label{c2}
 From the estimate (\ref{e10}) and Schwarz inequality, for any $\phi\in H^*$
 there exists $h\in H$ such that 
\begin{equation} \label{e11}
 |\phi(Z_2(t))|\leq |(Z_1(t),h)| \leq |Z_1(t)|\,|h|, \quad t\geq 0\;. 
\end{equation}
\end{cor}

\section{Inequalities} \label{sc4}

This is worth to emphasize the contributors to the maximal 
inequalities and exponential tail estimates. Kotelenez \cite{Ko1,Ko2} and Tubaro
\cite{Tu} studied the case of contraction semigroups when $p=2$. 
In \cite{Ko3} a maximal inequality for an analytic semigroup is
derived while Chow and Menaldi \cite{ChMe} obtained exponential 
tail estimates for some diffusion processes in Hilbert spaces.
\begin{theor} \label{th1}
Assume that processes $Z_1(t),~Z_2(t)$ and operators $S(t),~t\geq 0,$ are as
above and the process $\psi(t)$ fulfills ($\psi 1$). 
Then for any $\phi\in H^*,~p>2$, there exists a constant 
$\widetilde{c}_p>0$ such that 
\begin{equation} \label{e13}
 \mathbb{E}\left(\sup_{t\leq T}(\phi(Z_2(t)))^p\right) \leq \widetilde{c}_p
 \left(\sup_{t\leq T} || S(t)||^p\right)\;T^{p/2-1}\; \mathbb{E}\left( \int_0^T
 ||\psi(s)||_2^p\;ds\right) \;.
\end{equation}
\end{theor}
\proof{
 For $\alpha$ such that $1/p < \alpha < 1/2$, from  (\ref{e4a}), (\ref{e5}) 
 and (\ref{e3}) we have 
$$ 
   |Z_1(t)| = |I_\alpha Y_\alpha (t)|,~~~t\in[0,T]\;. 
$$
Then, by H\"older's inequality, where $q=p/(p-1)$:
\begin{eqnarray} \label{e14} 
 |Z_1(t)| &\leq& \frac{1}{\Gamma(\alpha)}\left| \int_0^t (t-s)^{(\alpha-1)}
 S(t-s)\,Y_\alpha (s)\,ds \right| \\
 &\leq& \frac{1}{\Gamma(\alpha)}\left( \int_0^t (t-s)^{(\alpha-1)q}
 ||S(t-s)||^q\,ds \right)^{\frac{1}{q}} \left(
  \int_0^t |Y_\alpha (s)|^p \,ds\right)^{\frac{1}{p}}\;. \nonumber
\end{eqnarray}
Now, from (\ref{e11}) and (\ref{e14}), 
$$ \left(\sup_{t\leq T}(\phi(Z_2(t)))^p\right) \leq 
 c_p \left( \int_0^T T^{(\alpha-1)q} ||S(s)||^q ds \right)^{\frac{p}{q}}
  \left( \int_0^T |Y_\alpha (s)|^p ds\right)\;,
$$
where $c_p:=\frac{|y|^p}{(\Gamma(\alpha))^p}$.  
Then
\begin{equation} \label{e15}
 \mathbb{E}\left(\sup_{t\leq T}(\phi(Z_2(t)))^p\right) \leq 
 c_p \left( \int_0^T T^{(\alpha-1)q} ||S(t)||^q dt \right)^{\frac{p}{q}}
 \mathbb{E}\left( \int_0^T |Y_\alpha (s)|^p ds\right) \;.
\end{equation}

Now, we shall estimate the last term in (\ref{e15}). We shall use
the following property of the stochastic integral: there exists
a constant $c$ that
\begin{equation} \label{e16}
 \mathbb{E}\left( \left| \int_0^t \psi(s)\,dW(s) \right|^p \right) 
 \leq c\, \mathbb{E}\left( \int_0^t ||\psi(s)||_2^2 \,ds\right)^\frac{p}{2} 
 ,~\mbox{where}~ p>0,~t\in[0,T].
\end{equation}
From (\ref{e4a}) and (\ref{e15}):   
\begin{eqnarray} 
  \mathbb{E}\left( \int_0^t |Y_\alpha (s)|^p\,ds \right) & \!\leq \!&
  \frac{c}{(\Gamma(1-\alpha))^p}\; \mathbb{E}\left\{ \int_0^T \left( \int_0^s
  ||(s-\sigma)^{-\alpha} S(s-\sigma) \psi(\sigma)||_2^2 \,d\sigma 
  \right)^\frac{p}{2} ds\right\} \nonumber \\  
  && \hspace{8ex} \mbox{(writing~out~the~Hilbert-Schmidt~norm)} \nonumber\\
  &\! \leq \!&
  \frac{c}{(\Gamma(1-\alpha))^p}\; \mathbb{E}\left\{ \int_0^T \left( \int_0^s
  || (s-\sigma)^{-\alpha} S(s-\sigma) ||^2||\psi(\sigma)||_2^2 \,d\sigma
  \right)^\frac{p}{2} ds\right\} \nonumber \\ & \!=\! &
  \frac{c}{(\Gamma(1-\alpha))^p}\; \mathbb{E}\left(\int_0^T 
  [(f*g)(s)]^\frac{p}{2}\,ds\right)\;, \nonumber 
\end{eqnarray}
where $f(s):= ||s^{-\alpha}\,S(s)||^2$, ~$g(s):= ||\psi (s)||_2^2$.

From the Young's inequality 
\begin{equation} \label{e17}
 \mathbb{E}\left( \int_0^T |Y_\alpha (s)|^p\,ds \right) \leq 
 \frac{c}{(\Gamma(1-\alpha))^p} \left( \int_0^T s^{-2\alpha} ||S(s)||^2 \,ds
 \right)^\frac{p}{2}\;\mathbb{E}\left( \int_0^T ||\psi(s)||_2^p\,ds \right) \;.
\end{equation}

The inequalities (\ref{e15}) and (\ref{e17}) provide
\begin{eqnarray} 
 \mathbb{E}\left(\sup_{t\leq T}(\phi(Z_2(t)))^p\right) & \leq &
 \widetilde{c}_p \left( \int_0^T  T^{(\alpha-1)q} ||S(t)||^q\,dt 
 \right)^\frac{p}{q} \left( \int_0^T t^{-2\alpha} ||S(t)||^2 \,dt
 \right)^\frac{p}{2} \nonumber \\ & &\; \times \;
 \mathbb{E}\left( \int_0^T || \psi (s)||_2^p\,ds
 \right) \nonumber\\ & \leq &
 \widetilde{c}_p \left( \sup_{t\leq T} ||S(t)||^p \right) T^{p/2-1}
 \mathbb{E}\left(\int_0^T || \psi (s)||_2^p\,ds \right)\;, \nonumber
\end{eqnarray}
where 
$$ \widetilde{c}_p := \frac{c\,c_p}{(\Gamma(1-\alpha))^p}, \quad 
\mbox{and} \quad
  p(\alpha-1)+\frac{p}{q}+(-2\alpha+1)\frac{p}{2}=\frac{p}{2}-1 \;.
$$
\hfill $\blacksquare$  } 

\noindent{\bf Comment:} Because the operators $S(t),~t\geq 0$, 
do not form any semigroup we can not expect the
equality of the processes $Z_1(t)$ and $Z_2(t)$, $t\geq 0$. In other words, 
by using the factorization method we are not able to prove that the
process $Z_2(t)$ has continuous modification $Z_1(t)$, $t\geq 0$. 

We may formulate the inequality "symmetric"to (\ref{e13}).
\begin{theor}
 Assume that $Z_1(t),Z_2(t)$ and operators $S(t),t\geq 0$, are as above. Then for any 
 $\phi\in H^*$, for arbitrary $p\in (2,\frac{1}{\alpha})$ and $\alpha\in (0,\frac{1}{2})$,
 there exists a constant $\hat{c}_p$ that
 \begin{equation} \label{e18}
  \mathbb{E}\left( \sup_{t\leq T}(\phi (Z_2(t)))^p\right) \leq 
  \hat{c}_p\left( \int_0^T t^{-2\alpha} || S(t)||_2^2\,dt \right)^{\frac{p}{2}}
  \mathbb{E}\left( \int_0^T ||\psi(t)||^p dt \right) \;.
 \end{equation}
\end{theor}
\proof{
The proof of (\ref{e18}) is nearly the same like the proof of (\ref{e13}).  
The different is the writing out the Hilbert-Schmidt norm 
$||(s-\sigma)^{-\alpha}\,S(s-\sigma)\psi(\sigma)||_2$ only.
\hfill $\blacksquare$  } \\

\noindent{\bf Comment:} 
In our case the assumptions of the Theorem 2 mean that the resolvent 
operators $S(t),t\geq 0$, must be of Hilbert-Schmidt type. 
The question is: what kind of Volterra
equations admit resolvents fulfilling that assumption? 
Good candidates seem to be the
integrodifferential equations (\ref{e8a}) mentioned earlier 
because of the form of the
resolvents. We can see (e.g.\ \cite{Pr}), that the resolvent 
operators $S(t),t\geq 0$, 
of (\ref{e8a}) are represented by the fundamental solutions 
$P_\alpha(t,x)$ of (\ref{e8a})
according to 
\begin{equation} \label{e19}
 (S(t)v)(x) = \int_{-\infty}^{\infty} P_\alpha(t,x-y)\,v(y)\,dy\;, 
 \quad t\geq 0,~x\in
 \mathbb{R}\;.
\end{equation}
Fundamental solutions $P_\alpha$ in (\ref{e19}) are well-known 
for $\alpha=1$, and for the
limiting cases $\alpha=0$ and $\alpha=2$. For our purposes, because
of the hypothesis (s), we may consider cases $\alpha\in [1,2)$, 
studied in details by 
\cite{Fu} and \cite{ScWy}.

Now, we shall adapt the result obtained by Peszat \cite{Pe} for
the convolution $Z_2(t),~t\in [0,T]$, given by (\ref{e6}).
We introduce the following definition and assumptions.
\begin{defin} \label{d1}
 We say that process $\psi : \Omega \times[0,T]\rightarrow L(H,H)$ is
 {\rm point-predictable} if for all $g,h\in H$ the process $(\psi(t)g,h),\,t\geq 0$ 
 is predictable
 with respect to the filtration $(\mathcal{F}_t),\,t\geq 0$.
\end{defin}
\noindent
($\psi 2$)\hfill \parbox[t]{150mm}{ \it
Here we assume that the Hilbert space $U=H$ and that the process 
$\psi : \Omega \times
[0,T]\rightarrow L(H,H)$ is point-predictable.}\\

\noindent
($\kappa$) \hfill \parbox[t]{150mm}{\it There exist $\alpha_0\in(0,\frac{1}{2})$ and 
$p_0>1$ such that 
$$ \kappa_T := \left( \int_0^T t^{(\alpha_0-1)p_0} ||S(t)||^{p_0} dt
\right)^{1/p_0} < +\infty\;,\quad \mbox{for any~} T>0.
$$ }
\begin{theor} \label{t3}
Assume that the operators $S(t),t\geq 0$, are as above, the process 
$\psi(t),t\geq 0$, fulfills ($\psi 2$), and conditions (s), ($\kappa$) are
satisfied. Assume that there exists a constant $\eta<+\infty$ such that 
\begin{equation} \label{e20}
 \sup_{0\leq t\leq T}\; \int_0^t (t-s)^{-2\alpha_0} ||S(t-s)\psi(s)||_2^2\,ds
 \leq \eta\,, \quad \mbox{P -- a.~s.}
\end{equation} 
Then for all $\delta >0$ there exist constants $C,\kappa_T$ that
\begin{equation} \label{e21}
 P\left\{ \sup_{0\leq t\leq T}\; |\phi(Z_2(t))| \geq \delta \right\} \leq 
 C \exp\left\{-\frac{\delta^2}{\kappa_T^2\,\eta} \right\} \;.
\end{equation}  
\end{theor}
\vspace{2mm}
 
\proof{
In our case, because we prove only the estimate (\ref{e21}) 
but not continuity, the proof is simple.
First of all, we formulate the inequality (\ref{e9}) 
in the case when the process 
$\psi(t),t\geq 0$, satisfies condition ($\psi 2$). 
Basing on Lemma 3.3 in \cite{Pe}, we have 
\begin{equation} \label{e22}
 \int_0^T \mathbb{E}\, \exp \left\{ \frac{1}{9\eta} |Y(t)|^2\right\}dt \leq 4T\;,
\end{equation}  
where 
$$ Y(t) := \int_0^t (t-s)^{-\alpha_0}S(t-s)\psi(s)dW(s), \quad t\in[0,T]\;.$$
Then, following the estimate from the proof of Lemma 3.4 in \cite{Pe}, we have:\\
for any $\phi\in H^*$ there exists $h\in H$ that 
\begin{equation} \label{e23}
 \sup_{0\leq t\leq T}\; \phi(Z_2(t)) \leq \frac{\sin \alpha_0\pi}{\pi}\,
 \sup_{0\leq t\leq T}\; |Z_1(t)||h|\leq \frac{c}{3}\,\kappa_T ||Y||_{L^q(0,T;H)}\;.
\end{equation} 
Now, using the estimates (\ref{e22}) and (\ref{e23}) we obtain the inequality 
$$ \mathbb{E} \exp \left\{\frac{\sup_{0\leq t\leq T}\; 
\phi(Z_2(t))}{\kappa_T^2\eta} \right\}\leq C\;. $$
This estimate and Doob's inequality complete the proof of Theorem \ref{t3}.
\hfill $\blacksquare$  } 

The factorization method applied to some class of stochastic linear Volterra
equations has provided similar estimates like that obtained earlier 
for stochastic evolution equations.
Unfortunately, because of lack of semigroup
property (in general, the resolvent operators $S(t),t\geq 0$, 
corresponding to the Volterra
equations considered do not form any semigroup), the factorization method does not 
provide continuity of stochastic convolutions arising in Volterra equations. 
The factorization method seemed to be more promising for obtaining continuity than it
finally appeared. We have seen, that 
the assumption about semigroup property is indispensable for
continuity by that method.

\end{document}